\documentclass[12pt]{article}
\usepackage{amsmath}
\newcommand{\be}{\begin{equation}}
      \newcommand{\ee}{\end{equation}}
      \newcommand{\ba}{\begin{eqnarray}}
       \newcommand{\ea}{\end{eqnarray}}
\newcommand{\ban}{\begin{eqnarray*}}
       \newcommand{\ean}{\end{eqnarray*}}

 \newcommand{\qed}{\hspace*{\fill}\rule{3mm}{3mm}\quad}
 \newcommand{\Pf}{\noindent {\em Proof.} }

\newcommand{\sect}[1]{\section{#1} \setcounter{equation}{0}}

\newtheorem{lem}{Lemma}[section]

\begin{document}
\newtheorem{defn}[lem]{Definition}
 \newtheorem{theo}[lem]{Theorem}
 \newtheorem{prop}[lem]{Proposition}
 \newtheorem{rk}[lem]{Remark}
 \newtheorem{ex}[lem]{Example}
 \newtheorem{note}[lem]{Note}
 \newtheorem{conj}[lem]{Conjecture}

\title{A Note On Obata's Rigidity Theorem I}
\author{Guoqiang Wu \& Rugang Ye \\ {\small Department  of Mathematics} \\
{\small University of California, Santa  Barbara} \\ \& 
\\
{\small Department of Mathematics, University of Science and Technology of China}} 
\date{March 2012}

\maketitle

Obata's rigidity theorem [O] as stated below is well-known.

\begin{theo} (Obata) Let $(M, g)$ be a connected complete
Riemannian manifold of dimension $n \ge 1$, which 
admits a nonconstant smooth solution of Obata's equation 
\ba \label{obata}
\nabla d w+wg=0.
\ea
Then $(M, g)$ is isometric to the $n$-dimensional round sphere ${\bf S}^n$.
\end{theo}

This theorem has various important geometric applications. For example, it is a main tool  for establishing the rigidity part of Lichnerowicz-Obata  
theorem [O1] [Be] regarding the first eigenvalue of the Laplacian under a positive lower bound for the Ricci curvature. 
Another example is that it leads to 
uniqueness of constant scalar curvature metrics in a conformal class of metrics containing Einstein metrics [O2][S]. On the other hand,  Obata's equation (\ref{obata}) stands out as an important and interesting geometric  equation for its own sake. 

In this note, we discuss various extensions of Obata's rigidity theorem. First we obtain general rigidity theorems and differentiable sphere theorems for the 
generalized Obata equation 
\ba \label{general}
\nabla d w+f(w)g=0
\ea
with a given smooth function $f$. Indeed, in Section 1 we first construct a class of Riemannian manifolds $M_{f, \mu}$ for a given function $f$ and a value 
$\mu$.  Then we derive a Jacobi field formula from the equation (\ref{general}), from which
the desired rigidity results with $M_{f, \mu}$ as 
model manifolds easily follow. The material in this section is adapted from the lecture notes [Y].  A  critical point of $w$ is assumed to exist in this section. In Section 2 
we formulate some natural conditions on $f$ and show that they imply existence of critical points of $w$.  

In Section 3 we present a derivation of global warped product 
structures implied by the equation (\ref{general}). Previous 
works on this subject have been done by Brinkmann in [Br] and 
Cheeger and Colding in [CC]. See that section for further discussions.  Materials in this section will be used in the subsequent sections and also in the sequel [WY] of this paper.

 In Section 4 we handle the hyperbolic case of the generalized Obata equation, i.~e.~the 
equation $\nabla dw-wg=0,$ and obtain hyperbolic versions of Obata's theorem.  In Section 5 similar results are obtained for 
the Euclidean cases of $(\ref{general})$.  
Our main rigidity results here involve  
a condition on the dimension of the solution space and improve Theorem 1.3 in [HPW] substantially.  Our results are actually optimal.  
In this context, besides analyzing the full rigidity   
case which requires the said dimension to be no less than a  
critical bound, we also characterize the geometry and 
topology of the underlying manifold when this dimension is lower. 

In the last section, we 
extend our results to the following more general formulation of (\ref{general})
\ba \label{more1}
\nabla d w+f(w, \cdot) g=0
\ea
with a given smooth $f$ defined on $I \times M$ for an interval 
$I$ and the general equation 
\ba \label{more2}
\nabla d w+zg=0,
\ea
where $w$ and $z$ are two smooth functions on $M$, which is equivalent to 
\ba \label{more3}
\nabla d w-\frac{\Delta w}{n} g=0,
\ea
where $n=dim\, M$. Another equivalent statement is that the Hessian of $w$ has only one eigenvalue everywhere. 
(For the general equation $\nabla d w= wq$ with an arbitrarily given smooth symmetric 2-tensor field $q$ and a smooth function $w$ we refer to [HPW], in which warped product rigidity is derived from this equation under a natural dimensional condition regarding its solution space.)

All manifolds in this note are assumed to be smooth. 
The  results in this note have been extended in 
[WY] to general Riemannian manifolds without the completeness assumption. 
Finally, we would like to acknowledge relavant discussions with T.~Colding,  G.~Wei and W.~Wylie.  

Note: We'll add additional references in the upcoming 
revised version of this paper. A further discussion of previous treatments of the equation (\ref{more3}) will also be provided.

%Previously, the following hyperblic and Euclidean extensions %of Obata's theorem 
%were obtained in [PW]. 

%\begin{theo} (Petersen and Wyllie) Let $(M, g)$ be a simply %connected 
%complete Riemannian manifold of dimension $n$. Assume that 
%the dimension of the space of smooth solutions of the %equation 
%\ba
%\nabla d w+ \tau w g=0
%\ea
%is at least $n+1$, where $\tau=-1$ or $0$. Then $(M, g)$ is 
%isometric to ${\bf H}^n$ if $\tau=-1$, and isometric to ${\bf %R}^n$ if 
%$\tau=0$.
%\end{theo}

%Our extension results inlude an improvement of this theorem.
%Namely we replace the lower bound $n+1$ by $n$ and remove 
%the simply connectedness condition. This improved result is 
%optimal. We also replace $n+1$ by $n-1$ and obtain the 
%corresponding regidity theorems. For details we refer 
%to Sections 2 and 3.

\sect{General Rigidity Theorems I}

 Before proceeding,  we would like to note that Obata's 
equation can be transformed to the equation $\nabla d w+
c wg=0$ for an arbitrary positive constant $c$ by a rescaling of the metric. This leads to an obvious rescaled extension  
of Obata's theorem. The same holds true for the various extensions of  Obata's equation in this paper.\\

\noindent {\bf 1.1  Rigidity Theorems} \\

\begin{theo} \label{general1} Let $(M, g)$ be a connected complete Riemannian manifold of dimension $n\ge 2$ which admits a nonconstant smooth solution $w$ of the 
generalized Obata equation (\ref{general})
for a smooth function $f(s)$. 
Assume that $w$ has at least one critical point 
$p$. Then $M$ is diffeomorphic to ${\bf  R}^n$ or ${\bf  S}^n$.
Moreover, 
$(M, g)$ is isometric to $M_{f, \mu}$ with $\mu=w(p)$.  
\end{theo}

The manifolds $M_{f, \mu}$  will be constructed below.  Theorem 
\ref{general1} leads to the following differentiable sphere theorem. 

\begin{theo} Let $(M, g)$ be an $n$-dimensional 
connected compact Riemannian manifold. Assume that it  admits a nonconstant smooth solution of the generalized 
Obata equation for some smooth function $f$.
 Then $M$ is diffeomorphic to ${\bf  S}^n$.
\end{theo}

Note that Obata's theorem easily follows from 
Theorem \ref{general1}, as it is easy to show that the
solution $w$ there must have a maximum point. \\

\noindent {\bf 1.2 Construction of $M_{f, \mu}$} \\

Let $f(s)$ be a smooth function defined on an interval $I=(A, B), [A, B), (A, B]$ or $(A, B)$ ($A$ is allowed to be $-\infty$ and $B$ is allowed to be $\infty$),  and let 
$\mu\in I$ satisfy $f(\mu)\not =0$.  
%(It is easy to extend our results to 
%suitable functions $f$ which are not necessarily defined on %the entire ${\bf R}$. For simplicity we present only the case %of $\bf R$.) 
Let $u$ be the unique maximally extended smooth solution of the initial 
value problem 
\ba \label{initial}
u''+f(u)=0, u(0)=\mu, u'(0)=0.
\ea 
By the uniqueness of $u$ we infer that 
$u$ is an even function. On the other hand, there holds $u''u'+f(u)u'=0$, which implies 
\ba \label{h-formula}
u'^2=-2h(u),
\ea
where $h$ is the antiderivative of $f$ such that $h(\mu)=0$.
%We present the construction in the case $f(\mu)<0$, while 
%that in the case $f(\mu)>0$ is similar. 
Let $T$ be the
supremum of $t$ such that $u$ is defined on $[0, t]$ and $u'
\not =0$ in $(0, t]$. We define for $n\ge 2$
\ba \label{metric}
g=dt^2+f(\mu)^{-2} u'^2 g_{{\bf  S}^{n-1}}
\ea
on $(0, T) \times {\bf  S}^{n-1}$. Set $\phi=-f(\mu)^{-1}u'$.
Then $\phi(0)=0, \phi'(0)=1$. Moreover, $\phi$ is an odd function because $u$ is even. It follows 
that $g$ extends to a smooth metric on 
the $n$-dimensional Euclidean open ball $B_T(0)$ where $B_T(0)-\{0\}$ is 
identified with $(0, T) \times {\bf  S}^{n-1}$. (If $T=\infty$, then 
$B_T(0)={\bf  R}^n$.)  There are three cases to consider.\\

\noindent {\it Case 1}   $h(s) \not =0$ for all $s>\mu$ if $f(\mu)<0$, and 
$h(s) \not =0$ for all $s<\mu$ if $f(\mu)>0$. \\
%. We shall say that $(f, \mu)$ is a {\it noncompact pair of %first kind}.\\

In this cas we make the following divergence assumption:
\ba
\int_{J} (-h)^{-1/2}ds=\infty,
\ea
where $J=(\mu, B)$ if $f(\mu)<0$ and 
$J=(A, \mu)$ if $f(\mu)>0$. 
To proceed, consider the function $s=u(t)$ 
on $[0, T)$. We have $t=u^{-1}(s)$. By (\ref{h-formula}) there
holds  
\ba
\frac{dt}{ds}=\pm (-2h(s))^{-1/2}
\ea
and hence
\ba \label{t-formula}
t=t(s)=\pm\int_{\mu}^s (-2h(\tau))^{-1/2} d\tau.
\ea
It follows that $T=\int_J (-2h)^{-1/2}ds=\infty$
and the manifold 
$({\bf  R}^n, g)$ is complete. We denote it by $M_{f, \mu}$. A pair $(f, \mu)$ satisfying the 
above conditions will be called {\it of noncompact type I}.\\

\noindent {\it Case 2} $h(\nu)=0$ for some $\nu$, where 
$\nu>\mu$ if $f(\mu)<0$ and $\nu<\mu$ if $f(\mu)>0$.  We assume that 
$\nu$ is the nearest such number from $\mu$. \\

\noindent {\it Case 2.1} There holds $f(\nu)=0$.  We say that $(f, \mu)$ is {\it of  noncompact type II}.\\

Since $h'(\nu)=f(\nu)=0$ it then follows that  
$\int_{J} (-h)^{-1/2} ds=\infty,$
where $J=(\mu, \nu)$ or $(\nu, \mu)$. 
As before, we infer that $T=\infty$. The complete 
manifold $({\bf  R}^n, g)$ is again denoted by $M_{f, \mu}$. \\

\noindent {\it Case 2.2} There holds $f(\nu) \not =0$. \\

\begin{lem} Assume $f(\nu) \not =0$. Then $T$ is contained in the domain of $u$ and $u(T)=\nu$. 
\end{lem}
\Pf We present the case $f(\mu)<0$, while the case 
$f(\mu)>0$ is similar. Since $h'=f$, the condition $f(\nu)\not =0$ implies 
$\int_{\mu}^{\nu} (-h)^{-1/2} ds<\infty.$  By the definition 
of $T$ there holds $u(t)<\nu$ for all $0\le t<T$. Hence 
 \ba
T \le \int_{\mu}^{\nu} (-2h(s))^{-1/2}ds<\infty.
\ea
On the other hand we have $|u''|\le |f(u)|\le \max \{|f(s)|: \mu\le s \le \nu\}$
on $[0, T)$. Consequently, $T$ is in the domain of $u$. By the definition of $T$ we then infer $T=\nu$. \qed \\

Next we assume in addition to $f(\nu) \not =0$ the {\it coincidence condition} $f(\mu)=-f(\nu)$. (The pair 
$(f, \mu)$ will be called {\it of compact type}.) Then 
the metric $g$ smoothly extends to ${\bf  S}^n$, where 
${\bf S}^n-\{p, -p\}$ ($p\in {\bf S}^n$) is identified with $(0, T) \times 
{\bf S}^{n-1}$. The Riemannian manifold $({\bf  S}^n, g)$ 
is also denoted by $M_{f, \mu}$. 

Note that the above arguments also provide a formula for the 
solution $u$. Indeed,  $u(t)$ is the 
inverse of the function $t(s)$ given by (\ref{t-formula}).

\begin{lem} In the above construction of $M_{f, \mu}$, the 
function $w=u(t)$ on $(0, T) \times {\bf S}^{n-1}$ smoothly 
extends to $M_{f, \mu}$ and satisfies the generalized 
Obata equation (\ref{general}) with the given $f$. 
\end{lem}
\Pf The evenness of $u$ implies that $u$ is a smooth function of $t^2$. But $t$ is the distance to the origin.  In
geodesic coordinates $x^i$ there holds $t^2=|x|^2$ and hence $w=u(t)$ extends smoothly across the origin. The situation in a second critical point is similar.  
On the other hand, 
a calculation similar to the proof of Lemma \ref{inverselemma} below shows that $w$ satisfies the generalized Obata equation on $(0, T) \times {\bf S}^{n-1}$, 
and hence on $M_{f, \mu}$. \qed \\

\noindent {\bf Examples 1} In the following examples the domain of $f$ is $\bf R$. \\
1) $f(s)=s^{2m}$ for a natural number $m$, $h(s)=(2m+1)^{-1} (s^{2m+1}-1)$ and $ 
\mu=1$. This is of noncompact type I.  \\
2) $f(s)=1, h(s)=s-1$ and $\mu=1$. This is also of noncompact type I. Note that 
$M_{1, 1}$ is the Euclidean space ${\bf R}^n$. \\
3) $f(s)=s^3-s, h(s)=\frac{1}{4}s^4-\frac{1}{2} s^2,
\mu=-\sqrt{2}$ and $\nu=0$. This is of noncompact type II. \\
4) $f(s)=s^{2m-1}$ for a natural 
number $m$, $h(s)=(2m)^{-1} (s^{2m}-1), \mu=1$ and $ \nu=-1$.  This is of compact type. 
Note that $M_{s, 1}$ is the round sphere $\bf S^n$. Indeed, 
there holds in this case $u=\cos t, u'=-\sin t$ and hence 
$g=dt^2+\sin^2 t \cdot g_{{\bf S}^{n-1}}$.   \\
5) $f(s)=\cos s, h(s)=\sin s, \mu=0, \nu=\pi$. This is of compact type. \\

\noindent {\bf 1.3 Calculation of Jacobi fields} \\

%For a point $p$ in a Riemannian manifold let $\Omega_p^T$ %denote the maximal star-shaped domain centered at the origin 
%on which $exp_p$ is defined. We set %$\Omega_p=exp_p(\Omega_p^T)$. 

\begin{lem} \label{gradlemma}
Let $w$ be a nonconstant smooth solution of the 
generalized Obata equation (\ref{general}) on a connected complete Riemannian manifold $(M, g)$ of dimension $n$ and 
for a smooth function $f(s)$. 
Let $p_0$ be a critical point of $w$.  Set $\mu=w(p_0)$. Then 
$f(\mu) \not =0$. Consequently, $p_0$ is a nondegenerate  local extremum point of $w$.  Moreover,  
if $\gamma$ is a unit speed geodesic starting at 
$p_0$, then there hold 
\ba \label{w-u}
w(\gamma(t))=u(t)
\ea
and 
\ba \label{gradient1}
\nabla w \circ \gamma=u' \gamma',
\ea
where $u$ is the solution of (\ref{initial}). In particular, 
$\gamma$ (with critical points of $w$ deleted) consists of  reparametrizations of gradient flow lines 
of $w$.
\end{lem}
\Pf Along each unit speed geodesic $\gamma(t)$ there holds 
\ba
\frac{d^2}{dt^2} w(\gamma(t))+f(w(t))=0.
\ea
We also have $\frac{d}{dt} w(\gamma(t))|_{t=0}=0$. Hence the formula (\ref{w-u}) holds true.
%Then there holds $w(p)=
%F(d(p_0, p))$ whenever there is a shortest geodesic from 
%$p_0$ to $p$. 
The claim $f(\mu) \not =0$ follows, because otherwise $w(\gamma(t)) \equiv \mu$ and hence $w\equiv \mu$ on $M$. 
The fact $f(\mu) \not =0$ and the equation (\ref{general}) imply that $p_0$ is 
a nondegenerate local extremum point of $w$. 

To see (\ref{gradient1}) we write $\nabla w \circ \gamma=X+\phi \gamma'$, where $X$ is normal. The generalized Obata equation is equivalent to 
\ba \label{nabla-equation}
\nabla_{\bf v} \nabla w+f(w){\bf v}=0
\ea
for all tangent vectors $\bf v$. We deduce that $X$ is parallel
and $\phi'+f(u)=0$. Hence we have  $X\equiv 0$ and 
$\phi'-u''=0$. The last equation and the initial values of $\phi$ and $u'$ imply $\phi=u'$. \qed \\
% Obviously, (\ref{gradient1}) implies  
%\ba \label{gradient2}
%\nabla w=u'(r) \frac{\partial}{\partial r}
%\ea
%inside of the injectivity domain at $p_0$, where $r$ denotes %the distance from $p_0$.

Next we present a calculation of Jacobi fields. 

\begin{prop} \label{jacobiprop}  Assume the same as in Lemma \ref{gradlemma}. 
Let $Y$ be a normal Jacobi field along $\gamma$ such that 
$Y(0)=0$ and $V$ denote the parallel transport of $Y'(0)$ along 
$\gamma$. Then there holds 
\ba \label{jacobi1}
Y=-f(\mu)^{-1} u'V.
\ea
\end{prop}
\Pf  Set ${\bf u}=
\gamma'(0)$, 
%Arguing by pulling back in terms of $exp_p$ we derive from %(\ref{gradient2}) the 
%formula (\ref{gradient1}) on $[0, T)$.    
%the metric to a suitable neighborhhod of the radial 
%line segment $\{ t{\bf u}: 0\le t \le T\}$, we can assume that 
%$\gamma(t), 0\le t \le T$ is contained inside of the 
%injectivity domain at $p$. Then (\ref{grad}) holds true in a %neighborhood 
%of $\gamma|_{[0, T]}$.   
${\bf v}=Y'(0)$ and $\gamma(t, s)=exp_{p_0}(t({\bf u}+s {\bf v})).$ Then 
\ba
Y(t)=\frac{\partial \gamma}{\partial s}|_{s=0}.
\ea
By (\ref{gradient1}) we deduce
\ba
\nabla w(\gamma(t, s))=|{\bf u}+s{\bf v}|^{-1}
u'(t|{\bf u}+s{\bf v}|) \frac{\partial \gamma}{\partial t}.
\ea
Hence 
\ba
\nabla_Y \nabla w &=& \frac{\partial }{\partial s}(|{\bf u}+s{\bf v}|^{-1}|\cdot
u'(t|{\bf u}+s{\bf v}|))|_{s=0} \frac{\partial \gamma}{\partial t}
\nonumber \\
&&+|{\bf u}+s{\bf v}|^{-1}
u'(t|{\bf u}+s{\bf v}|)\nabla_{\frac{\partial}{\partial t}}
Y.
\ea
Since $\frac{\partial}{\partial s} |{\bf u}+s{\bf v}|
= |{\bf u}+s{\bf v}|^{-1} {\bf u} \cdot {\bf v}=0$, we infer
\ba
\nabla_Y \nabla w =u'\nabla_{\frac{\partial }{\partial t}}
Y.
\ea
It follows that 
\ba
u'\nabla_{\frac{\partial }{\partial t}}
Y+f(u)Y=0.
\ea
Setting $Y=\phi V$ we deduce $u'\phi'+f(u)\phi=0$, 
i.~e.~$u'\phi'-u''\phi=0$. It follows that 
$\phi=Cu'$ for a constant $C$. Since $\phi(0)=0$ 
and $\phi'(0)=1$ we derive $\phi=-f(\mu)^{-1}u'$.
\qed \\

%From we infer that if the first conjugate point
%$p_1=\gamma(t_1)$  
%along $\gamma$ exists, then it is precisely the second point
%on $\gamma$ at which $\nabla w=0$. (Note that 
%$t_1$ is the first positive time at which $F_A'(t_1)=0$.) 
%Hence we 
%can apply the above result with $p_0$ replaced by $p_1$.
%Then the formulas (\ref{}) and (\ref{}) are extended beyond %the first conjugate 
%point. By induction, they are extended to the entire $\gamma$. %\qed \\

\noindent {\bf 1.4 Proof of  Theorem \ref{general1}} \\

\noindent {\it Proof of Theorem \ref{general1}} We can assume $n\ge 2$.  The formula 
(\ref{jacobi1}) implies
\ba \label{exp}
d exp_{p_0}|_{t{\bf u}}({\bf v})= -\frac{f(\mu)^{-1}u'}{t}V(t)
\ea  
which leads to 
\ba \label{metric}
 exp_{p_0}^*g=dr^2+f(\mu)^{-2}u'^2 g_{{\bf S}^{n-1}}
\ea 
for all $r>0$.
By Lemma \ref{gradlemma}, each unit speed geodesic $\gamma$ starting at $p_0$ is a reparametrization of a gradient flow line of $w$ (before reaching a critical point of $w$). 
Hence they can't meet each other before reaching a critical point of 
$w$. Moreover, no $\gamma$ can intersect itself before 
reaching a critical point of $w$.  By (\ref{gradient1}), they reach a critical point precisely 
at the first positive zero of $u'$.  If $u'$ has no 
positive zero, then we conclude that $exp_{p_0}$ is 
a diffeomorphism from $T_{p_0}M$ onto $M$. Next let $T$ be 
the first positive zero of $u'$. Then $exp_{p_0}$ is 
a diffeomorphism from $B_{T}(0)$ onto 
$B_{T}(p_0)$ (both are open balls). By    
(\ref{exp}), $exp_{p_0}$ maps $\partial B_{T}(0)$ 
onto a critical point $p_1$.   Employing the exponential map $exp_{p_1}$ we then infer that $exp_{p_0}$ extends to a smooth diffeomorphism from ${\bf S}^n$ onto $M$. Finally,  from the construction of $M_{f, \mu}$ and the completeness of $g$ it is easy to see that $(f, \mu)$ is one of the types in that construction and $(M, g)$ is isometric to $M_{f, \mu}$.
\qed \\

\sect{General Rigidity Theorem II} 

%In other words, $\Omega_p$ is 
%the maximal domain of $M$ whose points can be connected 
%to $p$ by geodesics.   
%Let $\rho_p$ denote the supremum 
%of $r>0$ such that $exp_p$ is defined on the open ball %$B_r(o)$ of $T_pM$. 
The main purpose of this section is to present natural conditions on $f$ which allow us to remove the condition of critical points in Theorem \ref{general1} and e.~g.~imply a differential sphere theorem without assuming compactness of the manifold.
\\

\noindent {\bf Definition 1} Let $f(s)$ be a smooth 
function on an interval $I$.  \\
1) We say that $f$ is {\it coercive}, if the following holds true.
Let $h$ be 
a somewhere negative antiderivative $h$ of $f$.
Then $h$ has zeros. Moreover, there holds 
$f(\mu)^2+f(\nu)^2 \not =0$ if $h(\mu)=h(\nu)=0$ and $h<0$ on $(\mu, \nu)$.  
If the maximal zero $\mu$ exists and 
$h<0$ on $I \cap (\mu, \infty)$, or  if the minimal zero 
$\mu$ exists and $h<0$ on $I \cap (-\infty, \mu)$, we also assume $f(\mu) \not =0$.  (A special case is that $f(\mu) \not =0$ for each zero $\mu$ of $h$.)  \\
2) We say that $f$ is degenerately coercive, if it is coercive and the following holds true. Let $h$ be an antiderivative $h$ of $f$.  If 
$h(\mu)=h(\nu)=0$ and $h<0$ on $(\mu, \nu)$, then 
$f(\mu)f(\nu)=0$.  \\
3) We say that $f$ is nondegenerately coercive, if the following holds true. Let $h$ 
be a somewhere negative antiderivative $h$ of $f$.
Then $h^{-1}((-\infty, 0))$ is a disjoint union of bounded intervals whose endpoints are contained in the domain of $f$. Moreover,  there holds $f(\mu)\not =0$ for 
each such endpoint $\mu$. \\ 

\noindent {\bf Examples 2} In the following examples, the domain of the function is ${\bf R}$. The function $f(s)=s^{2m-1}$ for a natural number $m$ is nondegenerately coercive. The function 
$f(s)=\pm s^{2m}$ for a  natural number $m$  is 
degenerately coercive. The functions $f(s)=1$ and $f(s)=1+\frac{1}{2} \cos s$ are degenerately coercive. 
The function $f(s)=\cos s$ is nondegenerately coercive. 
The function $f(s)=s^3-s$ is degenerately coercive.  
It is easy to construct many more examples. \\

\begin{theo} \label{coer1} Let $(M, g)$ be a connected complete Riemannian 
manifold of dimension $n\ge 2$ which admits a nonconstant smooth solution of the generalized Obata equation 
(\ref{general}) for a coercive function $f$. Then $M=M_{f, \mu}$ for some $\mu$.
In particular, $M$ is diffeomorphic to ${\bf  S}^n$ or ${\bf  R}^n$. 
\end{theo}

\begin{theo} \label{coer2} Let $(M, g)$ be a connected complete Riemannian 
manifold of dimension $n$ which admits a nonconstant smooth solution of the generalized Obata equation 
(\ref{general}) for a degenerately coercive function $f$. Then $M$ is diffeomorphic to ${\bf R}^n$. Moreover, if $n\ge 2$, then $(M, g)$ is isometric to $M_{f, \mu}$ for some $\mu$.
\end{theo}

\begin{theo} \label{coer3}
Let $(M, g)$ be a connected complete Riemannian 
manifold of dimension $n$ which admits a nonconstant smooth solution of the generalized Obata equation 
(\ref{general}) for a nondegenrately coercive function $f$. Then $M$ is diffeomorphic to ${\bf S}^n$. Moreover, if 
$n\ge 2$, then $(M, g)$ is isometric to $M_{f, \mu}$ for some $\mu$.
\end{theo}

In the case of Obata's equation, we have $f(s)=s$ and 
hence $h(s)=\frac{1}{2}s^2-C$. It is easy to see that 
$f$ is nondegeneratly coercive. Hence Obata's rigidity theorem is included in Theorem \ref{coer3} as a special case.

\begin{lem} \label{criticallemma1} Let  
$u$ be a nonconstant smooth solution of the equation 
\ba
u''+f(u)=0
\ea
on ${\bf R}$ for some smooth function $f$. Then the following hold true. \\
1) $u$ is  
symmetric with respect to each of its critical points. \\
2) If $f$ is coercive, then $u$ has at least one critical point. \\
3) If $f$ is degenerately coercive, then $u$ has precisely 
one critical point. \\
4) If  $f$ is nondegenerately coercive and $u'(t_0)=0$ for some $t_0$,  then 
$u'(t_1)=0$ for some $t_1>t_0$.  (It follows that $u$ is a periodic function.)  
\end{lem}
\Pf 1) Let $t_0$ be a critical point of $u$. Then 
$u(t_0-t)=u(t_0+t)$ follows from the uniqueness of the solution of (\ref{initial}). \\
2) Let $f$ be coercive. Assume that $u$ has no critical 
point. Then  $u({\bf R})$ is an 
open interval $(\mu_1, \mu_2)$. By (\ref{h-formula}) there holds 
$h \not =0$ on $(\mu_1, \mu_2)$. There holds for the inverse 
$t=t(s)$ of $u(t)$ 
\ba
t=\pm \int_{c}^s (-2h(s))^{-1/2}ds+t_0
\ea
with $c=u(t_0)$ for some $t_0$ in the domain of $u$. 
Since $u$ is defined on ${\bf R}$ we deduce
\ba
\int_c^{\mu_2} (-2h(s))^{-1/2}ds=\infty, 
\int_{\mu_1}^c (-2h(s))^{-1/2}ds=\infty.
\ea
If follows that $h(\mu_i)=0$ whenever $\mu_i$ is finite, 
$i=1,2$. It is impossible for both $\mu_1$ and $\mu_2$
to be infinite, otherwise $u({\bf R})={\bf R}$ and 
then $h(u(t))=0$ for some $t$, and hence $u'(t)=0$. 
By the coercivity assumption we then deduce that 
$h(\mu_i)=0$ and $f(u_i)\not =0$ for some $i$, say $i=2$.
But then 
\ba
\int_c^{\mu_2} (-2h(s))^{-1/2}ds<\infty
\ea
as $h'(\mu_2)=f(\mu_2) \not =0$. This is a contradiction. The case $i=1$ is similar. \\
3) Let $f$ be degenerately coercive. Assume that $u$ has 
more than one critical points.  Since $u$ is nonconstant, we translate the argument to achieve the following: $u$ has critical points $0$ and $t_0>0$, such that
$u' \not =0$ on $(0, t_0)$. Since $u$ is 
nonconstant, there holds $f(\mu) \not =0$, where $\mu=u(0)$.
It follows that $f(\nu)=0$, where $\nu=u(t_0)$. It follows then that
\ba
t_0=\int_I (-2h(s))^{-1/2}ds=\infty,
\ea
which is a contradiction, where $I=(\mu, \nu)$ or 
$(\nu, \mu)$.  \\
4) Let $f$ be nondegenerately coercive and $t_0$ a critical point of $u$. Translating the argument we can assume $t_0=0$. Then $u$ is an even function.
Hence it suffices to find one more critical point of $u$. Set $\mu=u(0)$. There holds 
$h(\mu)=0$. By (\ref{h-formula}), $h$ is nonpositive on the interval $u({\bf R})$. By the nondegenerately coercive assumption we then infer  that $u({\bf R}) \subset [a, b]$ for some finite $a$ and $b$ 
such that $h<0$ on $(a, b)$ and $h(a)=h(b)=0$. Obviously, there holds 
$\mu=a$ or $b$. We consider the former case, while the latter is similar. Assume that $0$ is the only critical point of 
$u$. Then $u({\bf R})=u([0, \infty))=[\mu, c)$ for some $c\le b$. It follows that
\ba \label{infinite}
\int_{\mu}^c (-2h(s))^{-1/2}ds=\infty.
\ea 
If 
$c<b$, then there holds $|h|>\delta$ in a neighborhood of 
$b$ for some 
$\delta>0$. Consequently we infer  
\ba \label{finite}
\int_{\mu}^c (-2h(s))^{-1/2}ds<\infty,
\ea
contradicting (\ref{infinite}).  
If $c=b$, we also derive (\ref{finite}) as $h'(b)=f(b)
\not =0$.  
\qed \\

\begin{lem}  \label{criticallemma2} Let $w$ be a nonconstant smooth solution of 
the generalized Obata equation with some $f$ on a complete Riemannian 
manifold $(M, g)$. 
If $f$ is coercive, then $w$ has at least one critical point. 
\end{lem}
%2) If $f$ is degenerately coercive, then  $w$ has at most one %critical point along each geodesic. \\
%3) Assume that $f$ is nondegenerately coercive and $p$ a %critical point of $w$. Then there is a $t_0>0$ with the %following properties.  Let $\gamma(t), t\ge 0$ be a geodesic %such that $\gamma(0)=p$. Then $\gamma(t_0)$ 
%is a critical point of $w$. Moreover, $\gamma(t)$ is not a %critical point of $w$ if $0<t<t_0$.  
%\end{lem}  
\Pf Choose a point $p \in M$ with $\nabla w(p) \not =0$. Let $\gamma$ be the  unit speed geodesic such that $\gamma(0)=
p$ and $\gamma'(0)$ is the unit vector in the direction of 
$\nabla w(p)$.  Set $u(t)=w(\gamma(t))$. Then there holds 
$u''+f(u)=0$. By Lemma, $u$ has at least one 
critical point $t_0$. Following the arguments in the proof of Lemma 
we deduce $\nabla w(\gamma(t))=\phi \gamma'$ 
with $\phi=u'+c$. But $u'(0)=|\nabla w(p)|
=\phi(0)$. Hence $\phi=u'$. If follows that $\nabla w(\gamma(t_0)=0$.
\qed \\

\noindent {\it Proof of Theorem \ref{coer1}}  This follows from 
Lemma \ref{criticallemma2} and Theorem \ref{general1}. \qed \\

\noindent {\it Proof of Theorem \ref{coer1} and Theorem \ref{coer2}} They follow from Lemma \ref{criticallemma2}, Lemma 
\ref{criticallemma1} and the proof of Theorem \ref{general1}.
In particular, the case of the exceptional dimension $n=1$
in Theorem \ref{coer1} follows from 2) of Lemma  
\ref{criticallemma1}, because a solution of (\ref{general1}) 
on ${\bf S}^1$ leads to a periodic function $u(t)$. 
\qed \\

\sect{Warped Product Structures} 

In [Br], using a result of partial differential equations and calculation in coordinates, Brinkmann derived from the equation (\ref{general}) a local warped product structure for 
the metric. In [CC], Cheeger and Colding derived from the equation (\ref{more2}) a global 
warped product structure for the metric in terms of  
calculation of differential forms. In this section, we 
present a slightly different derivation of  global warped product structures
based on the equation (\ref{general}).  Our approach 
uses the given solution $w$ as a global coordinate, which is motivated by the arguments in [Br].  This leads us to using the 
flow of the vector field $|\nabla w|^{-2} \nabla w$ in the construction. In comparison, the arguments in [CC] implicitly involve the vector field $|\nabla w|^{-1} \nabla w$. (The latter vector field also enters into our argument in an auxiliary and  different way, see the proof of Lemma \ref{h-equation1}.)
%One point of attention here is to construct a relevant global %product diffeomorphism on a given manifold (minus one or two %points), which is not  discussed 
%explicitly in [CC]. 
The results in this section will be extended to the general equation (\ref{more2}) in Section 6, where a special  
technical point regarding it will be handled.  (Some lemmas in this section are formulated for the general equation (\ref{more2}).)  
 
The formulation in this section  is particularly convenient for the applications in the subsequent sections. Its detailed arguments  are also needed in [WY] for dealing with incomplete manifolds.

%To understand why the critical point condition is necessary, %we 
%provide a characterization of manifolds admitting a %nontrivial 
%solution of Obata's equation. 

\begin{lem}  \label{level-set} Let $w$ and $z$ be two smooth functions on a Riemannian manifold $(M, g)$ satisfying the equation (\ref{more2}).
Let $N$  be a connected component of a level set of $w$.
Assume that $N$ contains no critical point of $w$.
Then $|\nabla w|$ and $z$ are constants on $N$. 
Moreover, the shape operator  of $N$ (with 
the normal direction given by $\nabla w$) is given by
$|\nabla w|^{-1} z_{N} Id$, where $z_{N}$ is the constant value of $z$ on $N$. In particular, $N$ is totally geodesic precisely when $z_{N}=0$.
\end{lem}
\Pf The equation (\ref{more2}) is equivalent to 
\ba \label{wz} 
\nabla_{\bf u} \nabla w+ z {\bf u}=0
\ea
 for all tangent vectors $\bf u$. We infer for $\bf u$ tangent to $N$
\ba
\nabla_{\bf u} |\nabla w|^2=2 \nabla_{\bf u} \nabla w \cdot \nabla w
=-2z {\bf u} \cdot \nabla w=0.
\ea
Hence $|\nabla w|$ is a constant on $N$. Next we have
\ba
\nabla_{\frac{\nabla w}{|\nabla w|^2}}
 \nabla w+ z \frac{\nabla w}{|\nabla w|^2}=0
\ea
and hence 
\ba
z=-\frac{1}{2} \nabla_{\frac{\nabla w}{|\nabla w|^2}} |\nabla w|^2.
\ea
Let $F(t, p)$ denote the flow lines of $\nabla w/|\nabla w|^2$ starting on $N$, where $p$ is an initial point and $t$ the time, with the initial value of $t$ being $\mu$, the value of $w$ on $N$.  There holds 
\ba
\frac{d}{dt} w(F(t, p))= \nabla w \cdot \frac{\nabla w}{|\nabla w|^2}=1.
\ea
Hence $w(F(t, p))=t$. By the established fact that $|\nabla w|$ is a constant on 
each component of any level set of $w$ we infer that 
$|\nabla w|(F(t, p))$ is independent of $p\in N$.  Now  
$\nabla_{\frac{\nabla w}{|\nabla w|^2}} |\nabla w|^2$ is equal to $\frac{d}{dt} |\nabla w|^2(F(t, p))|_{t=\mu}$, and hence independent of $p \in N$. We deduce 
that $z$ is a constant on $N$.
 
Finally we have 
\ba
\nabla_{\bf u} \frac{\nabla w}{|\nabla w|}
=-|\nabla w|^{-1} z {\bf u}
\ea
for tangent vectors $\bf u$ of $N$. 
\qed \\

%\noindent {\bf Definition 3} A tangent vector field $X$ with %isolated zeros
%on a Riemannian manifold is said to be {\it coercive}, if the %following holds true.  Let $\gamma(t)$ be a maximal flow line %of $|X|^{-1}X$ in the domain $X \not =0$. Then $\lim_{t %\rightarrow t_0} |X|
%=0$, where $t_0$ is a finite endpoint of $I$. \\

\begin{lem} \label{flow-line} Let $w$ and $z$ be two smooth functions on a Riemannian manifold $(M, g)$ satisfying the equation 
(\ref{more2}). 
Then the flow lines of $\nabla w/|\nabla w|$ in the domain $\{\nabla w \not =0\}$ are unit-speed
geodesics.  
\end{lem}
\Pf  By (\ref{wz}) we deduce 
\ba
\nabla_{\nabla w} \frac{\nabla w}{|\nabla w|}
&=& -|\nabla w|^{-1} z \nabla w-|\nabla w|^{-3} 
(\nabla_{\nabla w} \nabla w \cdot \nabla w) \nabla w
\nonumber \\
&=& -|\nabla w|^{-1} z \nabla w+|\nabla w|^{-1} 
z \nabla w=0.
\ea
The claim of the lemma follows. \qed \\

\begin{lem} \label{h-equation1} Assume that $w$ is a nonconstant solution of the 
generalized Obata equation (\ref{general}) on a Riemannian manifold $(M, g)$ for a given smooth function $f$.  Let $\mu$ be a value of $w$ and $N$ a connected component of $w^{-1}(\mu)$. 
Let $\alpha$ denote the value of $|\nabla w|$ on $N$, which is a constant by Lemma \ref{level-set}.   Then there holds for $p\in M$
\ba \label{hvalue}
|\nabla w|^2(p)=h(w(p)),
\ea
where $h(s)=\alpha^2-2\int_{\mu}^s f(\tau)d\tau$, as long as there is a gradient flow line $\gamma$ of $w$ such that $\gamma(t)$ converges to a point of $N$ in one direction and it converges to $p$ in the other direction. 
\end{lem} 
\Pf By a reparametrization we can assume that $\gamma$ is a flow line of 
$ \nabla w/|\nabla w|$, and hence a unit-speed geodesic by Lemma \ref{flow-line}. Set $u=w(\gamma(t))$. Then $u'^2=|\nabla w|^2$. 
As in Section 1, we have $u''+f(u)=0$ and  hence
$((u')^2-h(u))'=0$. (Note that $h$ here is  different from $h$ in Section 1.) We infer $(u')^2=h(u)+C$ or $|\nabla w|^2=h(w)+C$. Evaluating at a point of $N$ we deduce 
$C=0$. \qed \\

\begin{lem} \label{warplemma} Let $w$ be a nonconstant smooth solution of the generalized Obata equation (\ref{more2}) on 
a Riemannian manifold $(M, g)$ for a given smooth function $f$. Let $N$ be a connected component of $w^{-1}(\mu)$ for some $\mu$.
Assume that $\nabla w \not =0$ on $N$.  
As above, let $F(s, p)$ 
be the flow lines of $\nabla w/|\nabla w|^{2}$ starting on $N$ with the initial time being $\mu$. Assume that $F$ is smoothly defined on 
$I \times N$ for a time interval $I$.  Then there holds 
\ba \label{warp}
F^*g=\frac{ds^2}{h(s)}+\frac{h(s)}{\alpha^2} g_{N}
=\frac{ds^2}{h(s)}+\frac{h(s)}{h(\mu)} g_{N}, 
\ea
where $\alpha$ denotes the value of $|\nabla w|$ on $N$,  
$h$ is the same as in Lemma \ref{h-equation1} and $g_N$ is the induced metric on $N$. 
\end{lem}
\Pf Let ${\bf u} \in T_p N$ for some $p\in N$. Set 
$F_{\bf u}=dF_{(s, p)}({\bf u})$ and $F_s=\frac{\partial F}{\partial s}$. Using (\ref{nabla-equation}) and 
(\ref{hvalue}) we calculate 
\ba
\frac{\partial}{\partial s} |F_{\bf u}|^2
&=& 2\nabla_{F_{\bf u}}  F_s  
\cdot F_{\bf u}=2 \nabla_{F_{\bf u}} \frac{\nabla w}{|\nabla w|^2} \cdot F_{\bf u} \nonumber \\
&=& -2 f(w)|\nabla w|^{-2} |F_{\bf u}|^2=\frac{h'(s)}{h(s)}
|F_{\bf u}|^2.
\ea
Similarly, there holds
\ba
\frac{\partial}{\partial s} |F_s|^2
&=& 2\nabla_{\frac{\partial}{\partial s}} \frac{\nabla w}{|\nabla w|^2} \cdot F_s=-2f(w)|\nabla w|^{-2} |F_s|^2+
4 f(w)|\nabla w|^{-2} |F_s|^2 \nonumber \\
&=& 2f(w) |\nabla w|^{-2} |F_s|^2=-\frac{h'(s)}{h(s)}
|F_s|^2.
\ea
Integrating then leads to 
\ba
|F_{\bf u}|^2(p, s) &=& \frac{h(s)}{\alpha^2}|F_{\bf u}|^2(p, \mu),
\nonumber \\
|F_{s}|^2(p, s) &=& \frac{\alpha^2}{h(s)}|F_s|^2(p, \mu)
=\frac{1}{h(s)}.
\ea
\qed\\

\begin{theo} \label{brinkmann-global} Let $(M, g)$ be a connected complete Riemannian manifold which admits a nonconstant smooth solution $w$ of the generalized Obata equation 
(\ref{general}) for some smooth $f$.  Let $I$ denote the interior of the image $I_w$ of $w$. Let $\mu \in I$. Set $N=w^{-1}(\mu)$
and $\Omega=w^{-1}(I)$.   
Then $(N, g_N)$ is connected and complete with the induced metric 
$g_N$ and there is a diffeomorphism $F:
I \times N
\rightarrow \Omega$ such that $w(F(s, p))=s$ for all $(s, p)$.  The pullback metric $F^*g$ is a warped product metric  given by the formula 
(\ref{warp}).  Furthermore, there holds $M=\bar \Omega$
and $\partial \Omega$ consists of at most two points. 
Each point is either a unique global maximum point or a 
unique global minimum 
point of $w$.  
%In addition, 
%if $\mu$ is an endpoint of $\bar I$, then $w^{-1}(\mu)$ %consists of a single global extremum point of $w$. 

 Conversely, if $(N, g_N)$ is a Riemannian manifold and 
 $h(s)$ a positive smooth function on an interval $I$, then the function 
$w=s$ on $I \times N$ satisfies the generalized Obata equation with $f=-\frac{1}{2}h'$, where $I \times N$ is equipped with 
the metric $h^{-1}ds^2+ h g_N$.  
\end{theo}

 We would like to remark that this theorem can be used to 
 replace some arguments in Sections 1 and 2. This can e.~g.~be seen from the proofs of Theorem \ref{hyper2} and \ref{hyper3} 
 below. But the approach adopted  in these two sections is more concise. \\   

An immediate consequence of Theorem \ref{brinkmann-global} 
and Theorem \ref{general1} is the following result.

\begin{theo} \label{classify} Let $(M, g)$ be a connected complete Riemannian manifold which admits a nonconstant smooth solution $w$ of the generalized Obata equation 
(\ref{general}) for some smooth $f$. Then either $(M, g)$ is isometric to $M_{f, \mu}$ for some $f$ and $\mu$, or 
isometric to a warped product ${\bf R} \times_{\phi} (N, g_N)$ for a 
complete connected Riemannian manifold $(N, g_N)$ and a positive smooth function $\phi$ on ${\bf R}$. 
In the former case, $M$ is diffeomorphic to ${\bf R}^n$ if 
$w$ has precisely one critical point, and it is diffeomorphic to ${\bf S}^n$ if $w$ has two critical points. 
\end{theo}

\noindent {\it Proof of Theorem \ref{brinkmann-global}} \, Given the above results, the main point here is to construct the diffeomorphism $F$ in details, which requires some care in the case of a noncompact
$M$ and a possibly noncompact $N$.  
Let $N_0$ be a nonempty connected component of $N$ and let $F(s, p)$
be the same flow lines  as given in Lemma \ref{warplemma}, with $N$ repalced by $N_0$.  By the completeness of $g$, the induced metric 
$g_{N_0}$ is complete.
The formula $w(F(s, p))=s$ follows from a simple integration along the flow lines. Let $J_p$ be the interval of values of $w$ along the maximally defined $F(s, p)$ for $p\in N_0$.  
By Lemma \ref{level-set}, $|\nabla w|$ at $F(s, p)$ depends only on $s$. This fact and the completeness of $g$  imply that $J=J_p$ is independent of $p$.
Let $\Omega$ be the image of $F(s, p)$ for $p \in N_0$ and $s\in J$. Then $F: J \times N_0 \rightarrow \Omega$ is a diffeomorphism.
 
Let $p\in \partial \Omega$.
Then we have $F(s_k, p_k) \rightarrow p$ for some 
$s_k\in J$ and $p_k \in N_0$. There holds $s_k \rightarrow 
s^* \equiv w(p)$. We can assume that $s_k$ is a monotone
sequence. Let $\sigma(s)$ denote the value of $|\nabla w|$ 
at $F(s, p)$.  If $\nabla w(p) \not =0$, then 
$F(s, p_k)$ is defined on an open interval $J'$ containing 
$s^*$ as long as $k$ is large enough. Fix such a $k_0$.
Consider the case $s^*>\mu$, while the case $s^*<\mu$ is similar.  The length of the curve 
$\gamma_k(s)=F(s, p_{k}), \mu \le s \le s^*$ 
is given by  $L=\int_{\mu}^{s^*} \sigma^{-1}$. This integral is finite because $\sigma(s)$ is smooth and positive on 
$[\mu, s^*]$.  It follows that $dist(p, p_k)\le L+1$ for 
$k$ large. By the completeness of 
$(M, g)$ a subsequence of $p_k$ converges to 
a point $q\in N_0$. There holds $F(s, q) \rightarrow 
p$ as $s \rightarrow s^*$. But $F(s, q)$ is not defined at $s^*$, otherwise we would have $p \in \Omega$. 
We infer that $p$ is a critical point of $w$, and hence 
a nondegenerate local extremum point of $w$ (Lemma \ref{gradlemma}). 
It follows that the level sets of $w$ around $p$ are 
connected spheres filling a ball. (This is also clear from the proof of Theorem \ref{general1}.)  Since $F(s, q)$ passes through them, we infer that a neighborhood $U$ of $p$ satisfies $U-\{p\} \subset \Omega$.  Note that the image of 
$N_0$ under $F(s, \cdot)$ for $s$ close to $s^*$ is one of the said spheres.

Obviously, the above conclusion implies that $\bar \Omega$ is open. Hence $\bar \Omega=M$. We also infer $M=\bar \Omega
=\Omega \cup S$, where $S$ consists of at most two critical points of $w$. If $p\in S$, then $w(p)$ is an endpoint of $J=I$.
Moreover, $w(p_1) \not = w(p_2)$ if $S$ contains two points $p_1$ and $p_2$. All these also imply $N_0=N$. Finally, the claimed warped product formula 
(\ref{warp}) follows from Lemma \ref{warplemma}.
\qed \\

\sect{Hyperbolic Versions} 

\noindent {\bf 3.1 Main Theorems} \\

In this section we consider the following hyperbolic case  of the generalized Obata
equation 
\ba \label{hyper}
\nabla d w-wg=0.
\ea 
  
We first have the  following immediate consequence of Theorem \ref{general1}.

\begin{theo} \label{hyper-critical} Let $(M, g)$ be a connected complete Riemannian 
manifold which admits a nonconstant solution of equation 
(\ref{hyper-equation}) with 
critical points. Then it is isometric 
to ${\bf  H}^n$. 
\end{theo}

Note however that the function $f(s)=-s$ in (\ref{hyper}) is not coercive. Indeed, we can take the  
antiderivative $h(s)=-\frac{1}{2}s^2$.   Then $0$ is the only zero of $h$ and $h(0)=f(0)=0$.   More to the point, the solution 
$u=\sinh t$ of the equation 
\ba
u''-u=0
\ea
has no critical point. Hence Theorem \ref{coer1} and 
Theorem \ref{general1} are not applicable here. 
We consider instead the dimension of the solution space for 
(\ref{hyper}). 
First we would like to mention the following recent 
result ([Theorem 1.3, HPW]).

\begin{theo} (He, Petersen and Wylie) Let $(M, g)$ be a simply connected 
complete Riemannian manifold of dimension $n$. Assume that 
the dimension of the space of smooth solutions of the equation 
\ba \label{hyper-equation}
\nabla d w+ \tau w g=0
\ea
is at least $n+1$, where $\tau=-1$ or $0$. Then $(M, g)$ is 
isometric to ${\bf H}^n$ if $\tau=-1$, and isometric to ${\bf R}^n$ if 
$\tau=0$.
\end{theo}

We obtain the following two theorems and their 
Euclidean analogs which improve this result substantially and are indeed optimal. \\

\noindent {\bf Definition 2} For a Riemannian manifold 
$(M, g)$ and a smooth function $f(s)$ let $W_f(M, g)$ denote 
the space of smooth solutions of the generalized 
Obata equation (\ref{general}) on $(M, g)$.  We 
set $W_{h}(M, g)=W_{-s}(M, g)$.\\  

\begin{theo} \label{hyper1} Let $(M, g)$ be a connected complete Riemannian 
manifold of dimension $n\ge 2$.  Set $W_h=W_h(M, g)$.  Then 
$dim\, W_h\ge n$ iff 
$(M, g)$ is isometric to ${\bf  H}^n$.  Consequently, 
if $dim\, W_h \ge n$, then $dim \, W_h=n+1$. 
\end{theo}

\begin{theo} \label{hyper2} Let $(M, g)$ be a connected complete Riemannian 
manifold of dimension $n\ge 2$.  Set $W_h=W_h(M, g)$. Then $dim\, W_h$
$=n-1$ iff $(M, g)$ has constant sectional curvature -1 and 
is diffeomorphic to ${\bf R}^{n-1} \times {\bf S}^1$ (equivalently, 
$\pi_1(M)={\bf Z}$). More precisely, $dim\, W_h=n-1$ 
iff $(M, g)$ is isometric to 
${\bf H}^{n-1}_{\cosh}({\bf S}^1(\rho))$ or 
${\bf H}^{n-2}_{\cosh}({\bf H}_{\exp}({\bf S}^1(\rho)))$ for some $\rho>0$. (The former contains a closed geodesic while the latter doesn't.)
\end{theo}

The following theorem characterizes lower dimensions of $W_h(M, g)$.

\begin{theo} \label{hyper3} Let $(M, g)$ be a connected complete Riemannian 
manifold of dimension $n$ and $1\le k \le n-1$. Then $dim 
W_h(M, g)\ge k$ iff 
$M$ is isometric to ${\bf H}^{k}_{\cosh}(N, g_N)$ or 
${\bf H}^{k-1}_{\cosh}({\bf H}_{\exp}(N, g_N))$ for a connected 
complete Riemannian manifold $(N, g_N)$ of dimension $n-k$.
\end{theo}

The definition of the manifolds involving $\cosh$ and $\exp$ in the above theorems is given below. 
\\

\noindent {\bf Definition 3} Consider a Riemannian manifold 
$(N, g_N)$.\\ 
%We define its hyperbolic warpings ${\bf H}(N, g_N)$ as %follows.  First assume that $(N, g_N)$ is not isometric to %$S^m$. Then ${\bf H}(N, g_N)$ is defined to be any one of the 
%following three Riemannian manifolds $(\tilde N, g_{\tilde %N})$: \\
1) The cosh warping ${\bf H}_{\cosh}(N, g_N)$  of 
$(N, g_N)$ is defined to be the warped-product ${\bf R} \times_{cosh} (N, g_N)$. More precisely, it is defined to be   $(\tilde N, g_{\tilde N})$, where $\tilde N={\bf R} \times N$ and $g_{\tilde N}=
dr^2+\cosh^2 r \cdot g_N$. \\
2) The exponential warping ${\bf H}_{\exp}(N, g)=
(\tilde N, g_{\tilde N})$ is defined by   
$\tilde N={\bf R} \times N, g_{\tilde N}=dr^2+e^{2r} g_N$. \\
3) ${\bf H}^k_{\phi}(N, g_N)$ denotes the $k$-fold iteration of 
the $\phi$ warping, where $\phi=\cosh$ or $\exp$. \\ 
%3) The sinh warping $\tilde N=(0, \infty) \times N, 
%g_{\tilde N}=dt^2+\sinh^2 t \cdot g_N$. 
%If $(N, g_N)$ is isometric to $S^m$, then  ${\bf H}(N, g_N)$  
%is defined in the same way as above, except that the %definition of the sinh  warping is changed to ${\bf %H}^{m+1}$, 
%which is the metric completion of the manifold  in 3).
%The cosh warping is denoted by ${\bf H}_{cosh}(N, g_N)$. The 
%$k$-fold cosh warping  is defined to be the $k$-times iterated 
%cosh warping and denoted by ${\bf H}_{cosh}^k(N, g_N)$.
% as follows. Assume first that $(N, g_N)$ is not 
%isometric to $S^m$. Then it is defined to be 
%$(N\times R, \tilde g)$, where $\tilde g=
%ds^2+ \phi^2(s) g_N$, where $\phi(s)=\cosh s, \sinh s$ or 
%$e^s$. Next assume that $(N, g_N)$ is isometric to 
%$S^m$.   Let $(N, g_N)$ be a Riemannian manifold of dimension %$m$. Set $\tilde N=N \times R^k$ and 
More explicitly, ${\bf H}_{\cosh}^k(N, g_N)=(\tilde N, g_{\tilde N})$ with $\tilde N={\bf R}^k \times N$ and  
\ba \label{k-metric}
g_{\tilde N} &=& dr_1^2+\cosh^2 r_1 dr_2^2 
+\cdot \cdot \cdot+\cosh^2 r_1 \cdot \cdot \cdot \cosh^2 r_{k-1}
dr_{k}^2 
\nonumber \\
&&+\cosh^2 r_1 \cdot \cdot \cdot 
\cosh^2 r_{k} \cdot g_N.
\ea
Note that ${\bf H}_{\cosh}^{k}({\bf H}^m)={\bf H}^{m+k}$, 
${\bf H}^{n-1}_{\cosh}({\bf R})={\bf H}^n$, and ${\bf H}_{\cosh}^{n-1}
({\bf S}^1(\rho))$ is hyperbolic and diffeomorphic to ${\bf R}^{n-1} 
\times {\bf S}^1$, where ${\bf S}^1(\rho)$ denotes the circle of radius $\rho$. The last manifold contains a closed geodesic of length $2\pi \rho$. Furthermore, ${\bf H}_{\cosh}(N, g_N)$ is hyperbolic if $(N, g_N)$ is hyperbolic, and ${\bf H}_{\exp}(N, g_N)$ is hyperbolic if $(N, g_N)$ is flat.\\
%The $k$-fold cosh-hyperbolic warping of $(N, g_N)$ 
%is defined to be $(\tilde  N, g_{\tilde N})$ and 
%denoted by ${\bf H}_c^k(N, g_N)$. 

%The hyperbolic warpings ${\bf H}(N, g_N)$ of $(N, g_N)$ is %defined as follows. Assume first that $(N, g_N)$ is not 
%isometric to $S^m$. Then it is defined to be 
%$(N\times R, \tilde g)$, where $\tilde g=
%ds^2+ \phi^2(s) g_N$, where $\phi(s)=\cosh s, \sinh s$ or 
%$e^s$. Next assume that $(N, g_N)$ is isometric to 
%$S^m$.   
%Note that ${\bf H}^{k-1}({\bf H}(N, g_N))$ is isometric to 
%${\bf H}^n$ if ${\bf H}(N, g_N)$ is the sinh warping and 
%$(N, g_N)$ happens to be isometric to ${\bf S}^{n-k}$.\\
% to $R^{k-1} \times N$ and 
%(identifying $M$ with $R^{k-1} \times N$)
%\ba
%g &=& ds_1^2+\cosh^2 s_1 ds_2^2 
%+\cosh^2 s_1 \cdot \cdot \cdot \cosh^2 s_{k-2}
%ds_{k-1}^2 
%\nonumber \\
%&&+ \cosh^2 s_1 \cdot \cdot \cdot 
%\cosh^2 s_{k-1} \cdot c^2g_N,
%\ea
%where $c$ is a positive constant. The converse also holds 
%true.
%\\

\noindent {\bf 3.  Proofs of Main Theorems}\\

%The formulas in the following lemma are well-known and follow %from an elementary calculation. 

\begin{lem} \label{connection} Let $(N, g_N)$ be a Riemannian manifold and $\phi$ a positive smooth function on an interval 
$I$. Define $g=dr^2+\phi^2(r) g_N$ on $I \times N$. 
For a vector field $X$ on $I \times N$ which is tangent 
to $N$ we write $X=\sum a^{\alpha} \epsilon_{\alpha}$, where 
$\epsilon_{\alpha}$ is a local orthonormal frame of 
$(N, g_N)$. Let $\nabla^N$ be the Levi-Civita connection of $g_N$.  Then there hold w.~r.~t.~$g$
\ba
\nabla_{\bf v} X=\nabla^N_{\bf v} X-(X \cdot {\bf v}) 
\phi' \phi^{-1} \frac{\partial}{\partial r}\, \, \& \,\, \nabla_{\frac{\partial}{\partial r}} X=\sum (\phi a^{\alpha})' \epsilon_{\alpha},
\ea
where $\bf v$ is tangent to $N$. Consequently, there holds
$\nabla_{\frac{\partial}{\partial r}} (\psi X)=
(\phi \psi)'X,$
if $X$ is independent of $r$. On the other hand, there hold
for $\bf v$ tangent to $N$
\ba \label{r-formula}
\nabla_{\frac{\partial}{\partial r}} 
\frac{\partial }{\partial r}=0\,\, 
\& \,\, \nabla_{\bf v} 
\frac{\partial }{\partial r}=\phi' \phi^{-1} {\bf v}.
\ea
It follows in particular that each function $w=\int \phi(r)$ satisfies 
the generalized Obata (\ref{general}) equation with  
$f(s)=-\phi'(r(s))$, where $r(s)$ is the inverse of 
$\phi'(r)$.
\end{lem}
\Pf These formulas are well-known and follow from easy 
calculations. Using them it is easy to derive 
$\nabla dw=\phi'(r) g$, hence the claim regarding 
$w$ follows. (This has already been observed  in [CC].)
\qed \\

\begin{lem} \label{inverselemma} Let $(N, g)$ be a Riemannian manifold 
and $w_0$ a smooth solution of the equation (\ref{hyper-equation}) on $(N, g)$. Then 
the functions $\sinh r$ and $\cosh r \cdot w_0$ are solutions of 
(\ref{hyper-equation}) on ${\bf H}_{\cosh}(N, g)$. 
Consequently, $dim\, W_h({\bf H}_{\cosh}(N, g))
=dim\, (N, g)+1$.  
\end{lem}
\Pf  Set $\phi(r)= \cosh r$. Consider the function $w=\cosh r \cdot w_0=\phi(r) w_0$. We have 
\ba
\nabla w=\phi'w_0 \frac{\partial}{\partial r}+\phi^{-1} \nabla^N w_0.
\ea
By Lemma \ref{connection} we then deduce
\ba
\nabla_{\frac{\partial}{\partial r}} \nabla w
=\phi'' w_0 \frac{\partial}{\partial r}
+(\phi \phi^{-1})' \nabla^N w_0=w \frac{\partial}{\partial r} 
\ea 
and for ${\bf v}$ tangent to $N$
\ba
\nabla_{\bf v} \nabla w &=& \phi'({\bf v} w_0) \frac{\partial}{\partial r}+ (\phi')^2 \phi^{-1} w_0 {\bf v}
+\phi^{-1} \nabla^N_{\bf v} \nabla^N w_0
-({\bf v}w_0) \phi' \frac{\partial}{\partial r}
\nonumber \\
&=& (\phi')^2 \phi^{-1} w_0 {\bf v}+\phi^{-1} w_0 {\bf v}=
w {\bf v}.
\ea
The first claim of the lemma follows. 
Using (\ref{r-formula}) we also deduce that 
$\sinh r$ is a solution of (\ref{hyper-equation}). 
(This also follows from Lemma \ref{connection}.)  
\qed \\

The following lemma is an immediate  consequence.

\begin{lem} The functions $\sinh r_1, \cosh r_1 \sinh  r_2, 
\cosh r_1 \cosh r_2 \sinh r_3, ...\cosh r_1 \cdot \cdot \cdot 
\cosh r_{k-1} \sinh r_k$ and  $\cosh r_1 \cdot \cdot \cdot 
\cosh r_k \cdot w_0$ (c.~f.~(\ref{k-metric})) on 
${\bf H}^k_{\cosh}(N, g_N) $ are solutions of the equation 
(\ref{hyper-equation}), where $w_0$ is an arbitrary solution of (\ref{hyper-equation}) on $(N, g_N)$.
\end{lem}

\noindent {\it Proof of Theorem \ref{hyper3}} The ``if" part is provided by Lemma \ref{inverselemma}. Now we prove the 
``only if" part. We set $W=W_h(M, g)$ and assume  and assume $dim\,\, W\ge k$. The  equation (\ref{hyper-equation}) is linear, hence 
 $W$ is a vector space. As in [HPW] we choose $p_0\in M$ and consider the evaluation  map 
$\Phi_{p_0}: W \rightarrow {\bf R} \times T_{p_0}M, \Phi_{p_0}(w)=
(w(p_0), \nabla w(p_0))$.  By Lemma \ref{gradlemma}, its kernel is trivial.
Hence it is injective, as observed in [HPW]. Let $T_{p_0}M$ 
stand for $\{0\} \times T_{p_0}M$.   
We have $dim(im\, \Phi+T_{p_0}M)\le n+1$. Set $dim (im\, \Phi \cap 
T_{p_0}M)=l$. Then $dim(im\, \Phi+T_{p_0}M)=l+(k-l)+(n-l)=k+n-l.$
It follows that 
\ba
k+n-l\le n+1,
\ea
whence $l\ge k-1$. 

\noindent 1) First assume $k\ge 2$. Then $l\ge 1$. Choose  $w_0\in W$ such that 
$\Phi(w_0)$ is nonzero and belongs to $im\, \Phi \cap T_{p_0}M$. Then $w_0(p_0)=0, 
\nabla w_0(p_0) \not =0$. Set $N=w_0^{-1}(0)$. We apply 
Theorem \ref{brinkmann-global}. There holds $\mu=0$. We  choose $w_0$ such that $\alpha=|\nabla w_0(p_0)|=1$. Then we have  
\ba
h(s)=1-2\int_0^s (-\tau) d\tau=1+s^2.
\ea
On the other hand, the formula (\ref{warp}) implies 
%\ba
%\nabla w_1=h(t) \frac{\partial}{\partial %t}=(1+t^2)\frac{\partial}{\partial t}
%\ea
$|\nabla w_0|^2=|\nabla s|^2=h(s)=1+s^2 \ge 1$. It follows that 
$I={\bf R}$ 
and $F$ is a diffeomorphism from ${\bf R} 
\times N$ onto $M$. Setting $r=\sinh^{-1} s$  we deduce 
\ba \label{cosh-metric}
 F^*g=dr^2+ \cosh^2 r   \cdot g_N.
 \ea
 Moreover, $\sinh^{-1}$ maps $\bf R$ diffeomorphically onto $\bf R$. It follows that $(M, g)$ is isometric to ${\bf H}_{\cosh}
 (N, g_N)$.  
 
 By Lemma \ref{level-set} (or (\ref{cosh-metric})), $N$ is totally geodesic.  Hence 
 the restriction of each function in $W$ to $N$ satisfies the 
 equation (\ref{hyper-equation}) on $(N, g_N)$.
 Set $E_0=\{(a, {\bf v}) \in {\bf R} \times 
 T_pM: {\bf v} \perp \nabla w_0(p)\}$ and $W_0=\{
 w: w\in \Phi_p^{-1}(E_0)\}$. Then $dim \, W_0=k-1$.
 Note that $\nabla w(p_0)$ is tangent to $N$ for each $
 w\in W_0$. Applying the injectivity of the evaluation map 
 for $(N, g_N)$ at $p_0$ we infer that the restriction to $N$ maps $W_0$ injectively into $W_h(N, g_N)$.  
 It follows that $dim\, W_h(N, g_N) \ge k-1$. 
 If $k-1\ge 2$, we can repeat the above argument. By 
 induction, we deduce that $(M, g)$ is 
 isometric to  ${\bf H}_{cosh}^{k-1}(N, g_N)$ for 
 a connected complete  manifold $(N, g_N)$ of dimension 
 $n-k+1$. Moreover, $dim\, W_h(N, g_N) \ge 1$. \\ 
 \noindent 2) Next we consider the case $k\ge 1$. Note that 
 the case in 1) is reduced to this case at the last stage.
 Choose a function $w_0\in W_h(M, g)$ and a point $p_0\in M$ such that 
 $(w(p_0), \nabla w(p_0)) \not =(0, 0)$.  Then $w$ is
 nonconstant, for otherwise it would be zero.  If $\nabla w(p)=0$, we can apply Theorem \ref{hyper-critical} to infer that 
 $(M, g)$ is isometric to ${\bf H}^n={\bf H}_{\cosh}({\bf H}^{n-1})$. If $w(p)=0$, we can apply the above argument in 1) to deduce that $(M, g)$ is isometric to ${\bf H}_{\cosh}(N, g_N)$
 for a connected complete $(N, g_N)$. Finally we consider the 
 case $\mu \equiv w(p) \not =0, \nabla w(p)\not =0$. As before, we choose 
 $w$ such that $|\nabla w(p)|=1$.  We again apply Propostion \ref{brinkmann-global}. There holds 
 \ba
 h(s)=1+2\int_{\mu}^s  \tau d\tau=s^2+1-\mu^2.
 \ea
 There holds $|\nabla w|^2\ge 1-\mu^2$. If $\mu^2<1$, then $F$ is a diffeomorphism from 
${\bf R} \times N$ onto $M$. Setting $\sigma=\sqrt{1-\mu^2}$ and $r=\sinh^{-1} (t/\sigma)$ we deduce 
 \ba \label{dr}
 F^* g=dr^2+ \sigma^2 \cosh^2 r \cdot g_N.
 \ea
 Replacing $g_N$ by $\sigma^2 g_N$ we then infer that $(M, g)$ is isometric to ${\bf H}_{\cosh}(N, g_N)$. 
 If $\mu=1$, we set $s=e^r$ as long as $s>0$ 
 and deduce 
 \ba
 F^*g=dr^2+e^{2r} g_N.
 \ea
 Since $({\bf R} \times N, dr^2+e^{2r} g_N)$ is connected and complete, we infer that $s>0$ everywhere and  
 $F$ coupled with the function $e^r$ maps ${\bf R} \times N$ diffeomorphically onto $M$. It follows that 
 $(M, g)$ is isometric to ${\bf H}_{\exp}(N, g_N)$. 
 If $\mu=-1$, we set $s=-e^{r}$ and arrive at the same conclusion.  
 If $\mu^2>1$, we set $\sigma^2=\mu^2-1$ and set 
 $r=\cosh^{-1}(s/\sigma)$ to deduce 
 \ba
 F^* g=dr^2+ \sigma^2\sinh^2 r \cdot g_N.
 \ea
 There holds $|\nabla w|^2=\sigma^2 \sinh^2 r$.
 We infer that, as $r \rightarrow 0$, each geodesic in the $r$ direction converges to a critical point of 
 $w$. By Theorem \ref{hyper1} we conclude that $(M, g)$ is isometric to 
 ${\bf H}^n={\bf H}_{\cosh} ({\bf H}^{n-1})$.

\noindent 3) Combining the above two cases we then arrive at the 
claim of the theorem. \qed \\

\noindent {\it Proof of Theorem \ref{hyper2}} 
By Theorem \ref{hyper3} we infer that $(M, g)$ is isometric to 
the manifold ${\bf H}^{n-2}_{\cosh} ({\bf H}_{\phi}(N, g_N))$ 
with $\phi=\cosh$ or $\exp$. Then $N$ is 1-dimensional, and 
hence is isometric to either ${\bf R}$ or ${\bf S}^1(\rho)$ for 
some $\rho>0$. Moreover, ${\bf H}_{\phi}(N, g_N)$ is hyperbolic 
in either case of $\phi$. It follows that  ${\bf H}^{n-2}_{\cosh} ({\bf H}_{\phi}(N, g_N))$ is hyperbolic. \qed \\

\noindent {\it Proof of Theorem \ref{hyper1}} We follow the arguments in 1) of the proof of Theorem \ref{hyper3} and 
deduce that $(M, g)$ is isometric to 
${\bf H}^{n-1}_{\cosh}(N, g_N)$, where $(N, g_N)$ is 
either ${\bf R}$ or ${\bf S}^1(\rho)$ for some 
$\rho>0$.  We also deduce that 
$dim\, W_h(N, g_N)\ge 1$.  But $W_h({\bf S}^1(\rho))=\{0\}$.
Indeed, each $w\in W_h({\bf S}^1(\rho))$ can be given by the 
formula $w(\theta)=A\cosh (\theta+\theta_0)$ for some 
$A$ and $\theta_0$. Then $A=0$ because $w$ must be $2\pi$-periodic.  It follows that $(M, g)$ is isometric 
to ${\bf H}^{n-1}_{\cosh}({\bf R})={\bf H}^n$. \qed \\

\sect{Euclidean Versions}  

In this section we consider the following Euclidean 
analog of Obata's  
equation 
\ba \label{euc}
\nabla d w=0.
\ea

Let $W_e(M, g)$ denote the space of smooth solutions of 
(\ref{euc}) on $(M, g)$.

\begin{theo} \label{euc1} Let $(M, g)$ be an n-dimensional  connected complete Riemannian 
manifold. Then $dim\, W_e(M, g)\ge n$ iff $(M, g)$ is isometric 
to ${\bf  R}^n$ or ${\bf R}^{n-1} \times {\bf S}^1(\rho)$ for some $\rho>0$.
\end{theo}

Next we characterize the general situation $dim\, W_e(M, g) \ge k, 2\le k \le n$. 
Note $dim\,  W_e(M, g) \ge 1$ because of the presence of nonzero constant 
solutions. 

\begin{theo} \label{euc2} Let $(M, g)$ be an $n$-dimensional connected complete Riemannian manifold. Then $dim\, W_e(M, g) \ge k$ for $ 2\le k \le n$ iff 
$(M, g)$ is isometric to ${\bf R}^{k-1} \times 
N$ where $N$ is a connected complete Riemannian manifold.
\end{theo}
\Pf We argue as in 1) of the proof of Theorem \ref{hyper3}. 
Choose $w$ and apply Theorem 
\ref{brinkmann-global} in the same way as there. Now 
the function $h$ is given by $h(s)=1$ and hence 
\ba
F^*g=ds^2+g_N.
\ea
It follows that $(M, g)$ is isometric to ${\bf R} \times N$ 
with the product metric $ds^2+g_N$.  We apply the restriction to  $N$ and induction as before to arrive at the desired 
conclusion.  \qed \\

\noindent {\it Proof of Theorem \ref{euc1}} This theorem is 
an immediate consequence of Theorem \ref{euc2}. \qed \\

Another Euclidean version of the Obata equation is the following one.
\ba \label{euc2}
\nabla dw+g=0.
\ea

The following theorem is a special case of Theorem \ref{coer2}.
(It should be known.)

\begin{theo} Let $(M, g)$ be a connected complete 
Riemannian manifold of dimension $n\ge 1$, which admits 
a smooth solution of the equation (\ref{euc2}). Then 
$(M, g)$ is isometric to the Euclidean space 
${\bf  R}^n$.
\end{theo}

%\Pf First we claim that $w$ has at least one critical point. 
%Choose a point $p_0$. If $\nabla w(p_0)=0$, we are done. 
%So we assume $\nabla w(p_0)\not =0$. Consider unit speed %geodedics
%$\gamma$ with $\gamma(0)=p_0$. By the equation we deduce 
%$w''+1=0$ with $w=w(\gamma(t))$. It follows that 
%$w=-\frac{1}{2}t^2+at+b$, where $a=\gamma'(0) \cdot 
%\nabla w(p_0)$ and $b=w(p_0)$. Its maximum $\frac{1}{2}a^2+b$ %is achieved at $\gamma(a)$. Then the maximum of $w$ on $M$ 
%is achieved 
%at a point of $\gamma$,  where $\gamma'(0)$ is the unit %tangent vector in the direction of $\nabla w(p_0)$. 
%This point is a desired critical point. 
%Now we choose a critical point $p_0$ and apply Lemma. 
%It follows that $Y=t V$. It follows that $exp_{p_0}$ is a %local isometry from $(T_{p_0}M, g_{p_0})$ onto $(M, g)$.
%We claim that each unit speed geodesic starting at $p_0$ is 
%a ray, i.~e.~each portion of it is shortest.  Let $\gamma(t), 
%0\le t \le T$ be such a geodesic which fails to be shortest.
%Let $\gamma_1(t), 0\le t \le T_1$ be a shortest geodesic from  
%$p_0$ to $p$, where $T_1<T$. Then $w(\gamma(T))=
%w(p_0)-\frac{1}{2}T_1^2$ and $w(\gamma(T_1))=
%w(p_0)-\frac{1}{2}T^2$. This is a contradiction as 
%$\gamma(T)=\gamma(T_1)$.   It follows that $exp_{p_0}$ is 
%an isometry. \qed 

\sect{Rigidity Theorems for the General Equations (\ref{more1}), (\ref{more2}) and (\ref{more3})}

%It is natural to consider the following more general %formulation of (\ref{general})
%\ba \label{more1}
%\nabla d w+f(w, \cdot) g=0
%\ea
%for a smooth function $f(s, p)$ defined on $I \times M$, where 
%$I$ is an interval and $M$ a given manifold. Still more %general is the following version
%\ba \label{more2}
%\nabla d w+zg=0
%\ea
%for two smooth functions $w$ and $z$ on a given manifold.
%As it turns out, these two equations reduce to (\ref{general}).
%To save space, we formulate a combined result. 

%It is shown in [CC] that the function $z$ in (\ref{more}) 
%is related  to $w$ as follows. 
%We'll give a different proof of this fact, and determine ??
%\\

As mentioned in Section 3, a warped product  analysis of the equation (\ref{more2}) has been presented in [CC]. The focus in this section is on a delicate aspect of this equation concerning the relation between $z$ and $w$ around critical points of $w$.  We first formulate a few lemmas. 

\begin{lem} \label{newgradlemma} Let $(M, g)$ be a Riemannian manifold and $w$ and $z$ two smooth functions on $M$ satisfying 
(\ref{more2}). Let $p_0$ be a critical point of $w$ and 
$\gamma$ a unit speed geodesic with $\gamma(0)=p_0$. Then 
there holds 
\ba
\nabla w \circ \gamma=\frac{d w(\gamma(t))}{dt} \gamma'.
\ea
Hence the parts of $\gamma$ where $w(\gamma(t))'\not =0$ 
are gradient flow lines of $w$.
\end{lem}
\Pf The argument in the proof of Lemma \ref{gradlemma} can easily be adapted to the general equation (\ref{more2}).
\qed \\

\begin{lem} \label{isolated} Assume the same as in the above lemma. Assume in addition that $dim\, M\ge 2$ and $M$ is connected.  If $w$ has two sufficiently close 
critical points,  then 
it is a constant function. Consequently, if $w$ is a nonconstant 
function, then its critical points are isolated.
\end{lem}
\Pf  Let $p_0$ be a critical point of $w$ and $B_r(p_0)$ a convex geodesic ball. Assume that there is another critical point 
$p_0$ of $w$ in $B_r(p_0)$. Let $\gamma$ be the geodesic passing through $p_0$ and $p_1$. For $p\in B_r(p_0)-\gamma$,
let $\gamma_1$ be the shortest geodesic from $p_0$ to $p$, and 
$\gamma_2$ the shortest geodesic from $p_1$ to $p$.
Then they meet at $p$ nontangentially.  By Lemma \ref{newgradlemma},  we deduce that 
$p$ is a critical point of $w$. By continuity, every point in 
$B_r(p_0)$ is a critical point of $w$. The connectedness 
of $M$ then implies that every point of $M$ is a critical point of $w$, and hence $w$ is a constant. \qed \\

\begin{lem} \label{nondeg} Assume the same as in Lemma \ref{newgradlemma} and $dim\, M\ge 2$. Then each isolated critical point of $w$ is nondegenerate. Equivalently, $z(p_0) \not =0$ at each isolated 
critical point $p_0$ of $w$. 
\end{lem}
\Pf Let $p_0$ be an isolated critical point 
of $w$. 
By Lemma \ref{newgradlemma},  the geodesics starting at $p_0$
(with $p_0$ deleted) are 
reparametrizations of gradient flow lines of $w$ until 
they reach critical points of $w$. It follows that small geodesic 
spheres with center $p_0$ are perpendicular to 
the gradient of $w$ and hence are level sets of $w$. Consider a unit speed geodesic 
$\gamma(t)$ with $\gamma(0)=p_0$. Set $u(t)=w(\gamma(t))$. By the just derived fact and the connectedness of small geodesic spheres,  $u(t)$ is independent of the choice of $\gamma$.  By the arguments in the proof of 
Theorem \ref{general1}
we infer that the metric formula (\ref{metric}) holds true. If $u''(0)=0$, then 
the metric would be degenerate at $p_0$. 
Indeed, the $(n-1)$-dimensional volume of $\partial B_r(p_0)$ would be bounded from above by 
$cr^{2(n-1)}$. We conclude that $u''(0) \not =0$, which implies that $\nabla d w|_{p_0}$ is nonsingular.  \qed \\

%\begin{lem} Assume the same as in Lemma \ref{newgradlemma} and %$n\ge 2$, where $n$ is the dimension of $M$. Then there holds  %$z=f(w)$ for a smooth function $f$ defined on the image of $w$.
%\end{lem} 
%\Pf From the proof of Lemma we see that $u(t)$ is an even %function. Since $u''(0) \not =0$, we have %$u(t)=a+bt^2+O(t^4)$. Multiplying $w$ and $z$ by $-b^{-1}$ we %can assume $b=1$. Then 
%we have $b^{-1}(u-a)=t^2(1+O(t^2))$. By the implicit function %theorem we infer 
%$t^2=H(u-a)$ for a smooth function $H(s)$. Since $v$ is also %an even function, there holds $v=V(t^2)$ for a smooth function %$V$. 
%Hence $v=V(H(u-a))$. \qed \\

\begin{theo} \label{wztheorem} Let $w$ and $z$ be two smooth functions on 
a connected complete Riemannian manifold $(M, g)$ of dimension $n \ge 2$, such that  
(\ref{more2}) holds true. 
Then there is a unique smooth function $f$ on the image of $w$ such that $z=f(w)$. If  $w$  is nonconstant and 
has at least one critical point, then 
$M$ is diffeomorphic to either ${\bf R}^n$ or ${\bf S}^n$.
Moreover, $(M, g)$ is isometric to $M_{f, \mu}$ for some 
$\mu$, where $f$ is determined by the relation $z=f(w)$. 
If $w$ has no critical point, then $(M, g)$ is isometric to 
the warped product ${\bf R} \times_{\phi} (N, g_N)$ 
for a connected complete Riemannian manifold $(N, g_N)$ 
and a positive smooth function $\phi$ on ${\bf R}$.
\end{theo} 
\Pf The case of $w$ being a constant is trivial. So we assume that $w$ is a nonconstant function. Based on the above lemmas, it is clear that the proof of Theorem \ref{brinkmann-global} can be carried over to yield the same 
conclusions as there, without  the formula (\ref{warp}).
We set $f(s)=z(F(s, p))$ for a fixed $p\in N$, and 
$s\in I$, the interior of the image $I_w$ of $w$. Then $z=f(w)$ on $\Omega=F(I \times N)$.   Obviously, $f$ is smooth on $I$. 
It is also clear from the properties of $F$ that 
$f$ extends continuously to $I_w$. Let $\mu 
\in \partial I_w \cap I_w$. We claim that $f$ is smooth at $\mu$.
Let $p_0\in M$ be the unique critical point of $w$ such that 
$w(p_0)=\mu$. Let $\gamma$ be a unit speed geodesic with 
$\gamma(0)=p_0$ and set $u(t)=w(\gamma(t))$ as before. 
There holds $u(0)=\mu, u'(0)=0$. By Lemma \ref{nondeg}, $u''(0)\not =0$.
Since $u$ is independent of the choice of $\gamma$, it is an even function. Hence $u=\mu+at^2+Q(t^2)$ with $a \not =0$, 
where $Q$ is a smooth funtion with $Q(0)=Q'(0)=0$. By the implicit function theorem we deduce $t^2=H(u)$ for a smooth 
function $H$ and small $t$ and $u$.  On the other hand, $v(t)=z(\gamma(t))$ is also an even function because 
$z=f(w)$.
It follows that $v=G(t^2)$ for a smooth function $G$. We arrive at 
$v=G(H(u))$. Clearly, there holds $f(s)=G(H(s))$, and hence $f$ is smooth at $\mu$. 

With the above conclusion, the remaining part of the theorem follows from Theorem \ref{general1}. \qed \\  

\noindent {\bf Rmeark} In general, the above function $f$ is not smooth at critical values 
of $w$ if the dimension of $M$ is 1.  Consider e.~g.~$M={\bf R}, w=x^3$ and $z=-3x^2$.  Then $w$ and $z$ satisfy 
(\ref{more2}) and $f(s)=-3 s^{2/3}$. \\

An obvious equivalent formulation of the above result is the following theorem. 

\begin{theo} \label{Hessiantheorem} Let $w$ be a smooth function on 
a connected complete Riemannian manifold $(M, g)$ of dimension $n \ge 2$, such that  
(\ref{more3}) holds true. 
Then there is a unique smooth function $f$ on the image of $w$ such that $-\frac{1}{n}\Delta w=f(w)$. If  $w$  is nonconstant and 
has at least one critical point, then 
$M$ is diffeomorphic to either ${\bf R}^n$ or ${\bf S}^n$.
Moreover, $(M, g)$ is isometric to $M_{f, \mu}$ for some 
$\mu$, where $f$ is determined by the relation $-\frac{1}{n} \Delta w =f(w)$. If $w$ has no critical point, then $(M, g)$ is isometric to 
the warped product ${\bf R} \times_{\phi} (N, g_N)$ 
for a connected complete Riemannian manifold $(N, g_N)$ 
and a positive smooth function $\phi$ on ${\bf R}$.
\end{theo}

With the help of Theorem \ref{wztheorem}, all the results in 
Sections 2 and 3 concerning the generalized Obata equation 
(\ref{general}) extend to the more general equation 
(\ref{more1}). We  formulate this in a combined theorem as follows. 

\begin{theo} The theorems in Sections 2 and 3 continue to 
hold if the equation (\ref{general}) is replaced by 
the equation (\ref{more1}), and the conditions on $f(s)$
are assumed to hold for $f(s, p)$ for each fixed $p$. 
\end{theo}
 
To illustrate the more precise details, we also state an 
individual case explicitly as one example. 

\begin{theo} Let $(M, g)$ be a connected complete Riemannian manifold and $f$ a smooth function on $I \times M$ for an interval $I$. Assume that $f(\cdot, p)$ is nondegenerately 
coercive for each $p\in M$ and that there is a nonconstant solution of (\ref{more1}) on $(M, g)$. Then $M$ is diffeomorphic to ${\bf S}^n$. Moreover, if $dim\, M \ge 2$, then $(M, g)$ is isometric to $M_{f_0, \mu}$, where $f_0=f(\cdot, p_0)$ and $\mu=w(p_0)$ for a critical point $p_0$ of $w$.
\end{theo}

Finally, we state an easy consequence of the above results 
which provides a different angle of view. 

\begin{theo} Let $(M, g)$ be a connected complete Riemannian manifold and $f$ a smooth function on $I \times M$ for an interval $I$. Assume that there is a smooth solution $w$ of the 
equation (\ref{more1}). Then $f(s, \cdot)$ is a constant function on $M$ for 
each value $s$ of $w$. Consequently, no solution $w$ of (\ref{more1}) can exist for a generic 
$f$. 
\end{theo}

The condition of completeness can be removed, see [WY].

\end{document}